\input amstex
\documentstyle{amsppt}

\def\RC{{\text{\rm RatCurves}}^n(X)}

\def\F{\Cal F}
\def\O{\Cal O}
\def\C{\Cal C}
\def\U{\Cal U}
\def\K{\Cal K}
\def\V{\Cal V}
\def\Cx{\Cal C_x}

\def\P{\Bbb P^1}
\def\Tx{\Bbb P T_x(X)}
\def\spb{\smallpagebreak}
\def\mpb{\vskip 0.5truecm}
\def\bpb{\vskip 1truecm}

\def\spb{\smallpagebreak}
\def\mpb{\vskip 0.4truecm}
\def\bpb{\vskip 0.8truecm}

\magnification=\magstep1 
\NoRunningHeads \pagewidth{15 truecm} \pageheight{23 truecm}
\baselineskip=12pt
\parindent=20pt
\TagsOnRight
\def\spb{\smallpagebreak}
\def\mpb{\vskip 0.5truecm}
\def\bpb{\vskip 1truecm}
\vcorrection{-0.2 truecm} \hcorrection{0.5 truecm} \NoBlackBoxes
\topmatter
 \centerline{ \bf Birationality of the tangent map for minimal rational curves }
  \mpb \centerline{ \bf Jun-Muk Hwang \footnote{Supported by the Korea Research Foundation Grant
(KRF-2002-070-C00003).}  and Ngaiming Mok \footnote{Supported by a
CERG of the Research Grants Council of Hong Kong.}} \mpb
\centerline{\it Dedicated to Professor Yum-Tong Siu on his
sixtieth birthday}
\endtopmatter

\mpb

\centerline{Abstract} \medskip For a uniruled projective manifold,
we prove that a general rational curve of minimal degree through a
general point is uniquely determined by its tangent vector.  As
applications, among other things we give a new proof, using no Lie
theory, of our earlier result that a holomorphic map from a
rational homogeneous space of Picard number 1 onto a projective
manifold different from the projective space must be a
biholomorphic map.

\bpb
 \noindent {\bf \S 1. Introduction}

\spb Let $X$ be an irreducible uniruled projective variety. Let
$\RC$ be the normalized space of rational curves on $X$ in the
sense of [Ko]. For an irreducible component $\K$ of $\RC$, let
$\rho: \U \rightarrow \K$ and $\mu: \U \rightarrow X$ be the
associated universal family morphisms. In other words, $\rho$ is a
$\P$-bundle over $\K$ and for $\alpha \in \K$, the corresponding
rational curve in $X$ is $\mu(\rho^{-1}(\alpha))$. An irreducible
component $\K$ of $\RC$ is a {\it minimal component} if $\mu$ is
dominant and for a general point $x \in X$, $\mu^{-1}(x)$ is
projective. Members of a minimal component are called {\it minimal
rational curves}. For example, rational curves of minimal degree
passing through a very general point of $X$ are  minimal rational
curves. Here, `very general' means that the point is chosen
outside a countable union of subvarieties of dimension strictly
smaller than the dimension of $X$.  Denote by $\Bbb P T(X)$ the
projectivization of the tangent bundle of the smooth part of $X$.
Given a minimal component $\K$, consider the rational map
$$\tau: \U \dashrightarrow \Bbb P T(X)$$ defined by
$$\tau(\alpha) := \Bbb P T_x(C)$$ for $\alpha \in \U$ such that
$x:=\mu(\alpha)$ is a smooth point of $X$ and $\rho (\alpha)$
corresponds to a rational curve $C$ on $X$ smooth at $x$.  Let $\C
\subset \Bbb PT(X)$ be the proper image of $\tau$. We call $\tau$
the {\it tangent map} of $\K$ and $\C$ the {\it total variety of
minimal rational tangents} of $\K$.

When $X$ is smooth, it is well-known that the degree of a minimal
rational curve with respect to the anti-canonical bundle is
bounded by $\dim (X) +1$ and the tangent map $\tau: \U
\dashrightarrow \C$ is generically finite. We will prove the
following.

\spb \proclaim{Theorem 1} For any uniruled projective manifold $X$
and any minimal component $\K$, the tangent map $\tau: \U
\dashrightarrow \C$ is birational.
\endproclaim

\spb In the process of proving Theorem 1, we will also prove

\spb \proclaim{Theorem 2} Suppose a uniruled projective manifold
has two distinct minimal components $\K$ and $\K'$. Then their
total spaces of minimal rational tangents $\C$ and $\C'$ are
distinct. \endproclaim

\spb Theorem 1 and Theorem 2 imply that
  a general rational curve of minimal degree through a general point
of $X$ is uniquely determined by its tangent vector. Theorem 1 is
proved under some special assumptions on $\K$, which do not always
hold,  in [KK].  Their method is completely different from ours.

 In addition to its intrinsic interest, Theorem 1 has several
interesting consequences. For a  point $x \in X$, the subvariety
$\C_x = \C \subset \Bbb PT_x(X)$ is called the {\it variety of
minimal rational tangents} at $x$. It is known that $\mu^{-1}(x)$
is a (not necessarily irreducible) smooth projective variety for a
general point $x$.  Kebekus showed that for a general $x$, the
restriction of the tangent map to the fiber at $x$, $\tau_x:
\mu^{-1}(x) \rightarrow \Cx$ is a finite morphism ([Ke, Theorem
3.4]). Thus Theorem 1 has the following consequence.

\spb \proclaim{Corollary 1} For any uniruled projective manifold
and a minimal component, the normalization of the variety of
minimal rational tangents at a general point is smooth.
\endproclaim

\spb This has an application  to the rigidity of generically
finite morphisms to Fano manifolds of Picard number 1. We say that
$\C_x$ is non-linear, if one (hence all) of the components of
$\Cx$  is not a linear sub-variety of $\Tx$. We will prove

\spb \proclaim{Theorem 3}  Let $\chi: {\Cal X} \rightarrow
\Delta:= \{ t \in {\Bbb C}, |t|<1 \}$
 be a regular family of Fano
manifolds of Picard number 1 so that $X_{0}= \chi^{-1}(0)$ has a
minimal component whose variety of minimal rational tangents at a
general point is non-linear.
 For a given projective manifold $X'$,
 suppose there exists a surjective morphism
 $f: {\Cal X}'= X' \times \Delta \rightarrow {\Cal X}$ respecting
 the projections to $\Delta$ so that
$f_{t}: X' \rightarrow X_{t}$ is a  generically finite morphism
for each $t \in \Delta$. Then there exists $\epsilon >0$ and a
holomorphic family of biholomorphic morphisms $v_{t}:X_{o}
\rightarrow X_{t}$ for $|t| <\epsilon$,
 satisfying $v_{0}=\text{{\rm Id}} $ and $f_{t} = v_{t} \circ f_{o}$. \endproclaim

\spb One interesting consequence is the following.

\spb \proclaim{Theorem 4} Let $G/P$ be a rational homogeneous
space of Picard number 1. If $f:G/P \rightarrow X$ is a surjective
morphism to a smooth projective variety $X$ of positive dimension,
then either $X \cong \Bbb P^n, n = \dim(G/P)$, or $f$ is a
biregular morphism. \endproclaim

\spb This was originally proved in [HM1] by using quite a bit of
Lie theory. The proof given here uses  no Lie theory.

 A natural question arising from Theorem 2 is to understand when a projective manifold
  $X$ admits a
minimal component $\K$ whose variety of minimal rational tangents
at a general point is non-linear. We expect the following.

\spb \proclaim{Conjecture } For any Fano manifold of Picard number
1 excepting the projective space and for any minimal component,
the variety of minimal rational tangents  at a general point is
non-linear if its dimension is positive. In particular, for a Fano
manifold of Picard number 1 with index $>2$, any minimal component
has  non-linear varieties of minimal rational tangents at general
points.
\endproclaim

\spb There are two ingredients of the proof of Theorem 1, the work
of Cho-Miyaoka-Shepherd-Barron ([CMS]) and the theory of
differential systems on the subvariety of the projectivized
tangent bundle. Roughly speaking, [CMS] proves Theorem 1 in the
special case when $\tau: \U \dashrightarrow \Bbb PT(X)$ is
dominant and the theory of differential systems reduces the
general case to this special case. The theory of differential
systems we need will be explained in Section 2. Theorem 1 and
Theorem 2 will be proved in Section 3. Theorem 3 and Threorem 4
will be proved in Section 4. We will work over the field of
complex numbers.

\mpb
 \noindent {\bf \S 2. Differential
systems on the subvariety of the projectivized tangent bundle}

\spb We need to study some natural distributions defined on a
subvariety of the projectivized tangent bundle of a complex
manifold.  Since this theory does not seem to be well-known to
algebraic geometers or complex analysts, we will give a full
account in this section. In the standard literature, a
distribution on a complex manifold means a subbundle of the
tangent bundle.  In this article, it will be used in a broader
sense.

A
 {\it distribution} on an irreducible
normal variety is a subbundle of the tangent bundle of a Zariski
open subset of the smooth part of the normal variety. This Zariski
open subset will be called the {\it domain of definition} of the
distribution. The complement of the domain of definition is called
the {\it singular locus} of the distribution. For simplicity, we
will regard two distributions identical if they agree on the
intersection of their domains of definition.

 By taking local holomorphic sections, we can view a
distribution $\Cal V$ as a locally free sheaf on the domain of
definition. The Lie bracket of local vector fields define a $\Bbb
C$-linear sheaf map $[,]: \Cal V \times \Cal V \to \Theta$ where
$\Theta$ denotes the tangent sheaf of the domain of definition.
The distribution defined by $[\Cal V, \Cal V] + \Cal V$ is called
the {\it first derived system} of $\Cal V$ and is denoted by
$\partial \Cal V$. The domain of definition of $\partial \Cal V$
may be different from that of $\Cal V$. The Lie bracket defines a
holomorphic vector bundle morphism $\bigwedge^2 \Cal V \to
\Theta/\Cal V$ over the domain of definition. This section of
$Hom(\bigwedge^2 \Cal V, \Theta/\Cal V)$ over the domain of
definition will be called the {\it Frobenius bracket tensor} of
$\Cal V$. By abuse of notation, we will use the same notation
$[,]$ for the Frobenius bracket tensor and the Lie bracket of
vector fields. By the famous Frobenius theorem, if the Frobenius
bracket tensor vanishes, or equivalently, if $\partial \Cal V =
\Cal V$, the distribution $\Cal V$ is integrable and defines a
foliation on its domain of definition. A {\it leaf} of an
integrable distribution means a maximal submanifold  in the domain
of definition which is tangent to the given distribution. The {\it
Cauchy characteristic} of $\Cal V$ is the distribution $Ch(\Cal
V)$ which is defined at a general point $x$ by
$$Ch(\Cal V)_x = \{ v \in \Cal V_x, [v, \Cal V_x] = 0 \}. $$
The Cauchy characteristic of any distribution is integrable, as
can be easily checked  using the Jacobi identity.

 Let $\xi: Y \rightarrow X$ be a dominant morphism between two
varieties. For a distribution $\Cal D$ on $X$,  $d \xi^* \Cal D$
denotes the distribution on $Y$ defined by $$ (d \xi^* \Cal D)_y =
(d \xi_y)^{-1}( \Cal D_{\xi(y)})
$$ at a general point $y \in Y$ such that the differential $d
\xi_{y} : T_y(Y) \to T_{\xi(y)}(X)$ is surjective and $\xi(y)$
lies in the domain of definition of $\Cal D$. $d \xi^* \Cal D$
will be called {\it the pull-back of } $\Cal D$ by $\xi$.  Given a
local section $w$ of $\Cal D$ in a neighborhood of $\xi(y)$, we
can find a local section $d \xi^* w$ of $d \xi^* \Cal D$ in a
neighborhood of $y$ so that $w = d \xi_* (d \xi^* w)$. We will say
that $d \xi^* w$ is a {\it lift} of $w$. Two different lifts of
$w$ differ by a local vector field tangent to the fibers of $\xi$.
The proof of the following lemma is elementary.

\spb \proclaim{Lemma 1} For any dominant morphism $\xi:  Y \to X$
and a distribution $\Cal D$ on $X$, $$ \align  d \xi^*(
\partial \Cal D) &= \partial (d \xi^* \Cal D)\\   d \xi^*
Ch( \Cal D) &= Ch(d \xi^* \Cal D). \endalign$$ In particular,
tangent vectors to the fibers of $\xi$ are contained in $Ch(d
\xi^* \Cal D).$ \endproclaim

\spb For the rest of this section, let us fix a complex manifold
$X$ and
 a subvariety  $\Cal C \subset \Bbb P T(X)$ which is dominant over
 $X$. There are two naturally defined distributions ${\Cal J}$ and ${\Cal P}$ on ${\Cal
 C}$.
At a general point $\alpha \in \C$, they are defined  by
$$ \align
{\Cal J}_{\alpha} &:= (d \pi)^{-1}(\Bbb C \alpha) \\
{\Cal P}_{\alpha} &:= (d \pi)^{-1}( \hat T_{\alpha}(\C_x))
\endalign $$
where $d \pi: T_{\alpha}({\Cal C}) \rightarrow T_{x}(X)$ is the
differential of the natural projection $\pi: {\Cal C} \rightarrow
X$ at $\alpha \in {\Cal C}, x= \pi(\alpha)$, and $\hat
T_{\alpha}(\Cx) \subset T_{x}(X)$ is the linear tangent space of
${\Cal C}_{x}:= \pi^{-1}(x) \subset {\Bbb P}T_{x}(X)$ at $\alpha$.
Both ${\Cal J}$ and ${\Cal P}$ are canonically determined by
${\Cal C}$. ${\Cal J}$ has rank $p + 1$ and ${\Cal P}$ has rank
$2p+1$, where $p$ is the fiber dimension of $\pi: {\Cal C}
\rightarrow X$. Also we have the trivial vertical distribution
${\Cal V}$ of rank $p$ on ${\Cal C}$ defining the fibers of $\pi$.
Clearly, ${\Cal V} \subset {\Cal J} \subset {\Cal P}$.

\spb \proclaim{Proposition 1} In a neighborhood of a general point
of $\C$, choose a line subbundle $\Cal F$ of $\Cal J$ such that
$\Cal V + \Cal F = \Cal J$. Then a local section of $\Cal P$ is of
the form $[v , f ] + u + f'$ where $u,v$ are local sections of
$\Cal V$ and $f, f'$ are local sections of $\Cal F$. In other
words, $\Cal P =
\partial \Cal J$. \endproclaim

\spb {\it Proof}. For notational simplicity, we will work over
$\Xi: = T(X)\setminus($0-section$)$. Let $\xi: \Xi \rightarrow
\Bbb P T(X)$ be the natural $\Bbb C^{*}$-bundle. Let $\hat {\Cal
C} := \xi^{-1}({\Cal C})$. We will denote the restriction of $\xi$
to $\hat{\Cal C}$ by the same letter $\xi$. Let $\hat{\Cal J}:= d
\xi^{*}{\Cal J}, \hat{\Cal P}:= d \xi^{*}{\Cal P}, \hat{\Cal V} :=
d \xi^{*}{\Cal V}$, and $\hat{\Cal F} :=d  \xi^{*}{\Cal F}$ be the
distributions pulled-back to $\hat{\Cal C}$. By Lemma 1, it
suffices to check that $\hat{\Cal P} = \partial \hat{\Cal J}$.

We start with $ \partial \hat{\Cal J} \subset \hat{\Cal P}$. It
suffices to show $[\hat{\Cal V}, \hat{\Cal F}] \subset \hat{\Cal
P}$. Let $x_{1}, \ldots, x_{n}$ be a local coordinate system on
$X$. Let $\lambda_{1}=d x_{1}, \ldots, \lambda_{n}= d x_{n}$ be
linear coordinates in the vertical direction of $\Xi$. Let $\hat
v=\sum_{i}v_{i} \frac{\partial}{\partial \lambda_{i}}$ be a lift
of a  local section of ${\Cal V}$ and $\hat f=
\sum_{i}f_{i}\frac{\partial}{\partial \lambda_{i}} + \zeta
\sum_{j}\lambda_{j} \frac{\partial}{
\partial x_{j}}$ be a lift of a local section of ${\Cal F}$
over a small open set in $\hat{\Cal C}$. Here $v_{i}, f_{i},
\zeta$ are suitable local holomorphic functions. Dividing by
$\zeta$ and looking at generic points outside the zero set of
$\zeta$, we may assume that $\zeta \equiv 1$. Then $[\hat v, \hat
f] = \sum_{i}v_{i} \frac{\partial}{\partial x_{i}}$ modulo
$\hat{\Cal V}$. But this is precisely the vectors $\hat v$ viewed
as the tangent vectors to $X$. Hence $[\hat v, \hat f] \in
\hat{\Cal P}$.

{}From the above expression of $[\hat v, \hat f]$ modulo
$\hat{\Cal V}$, we see that the rank of $\partial \hat{\Cal J}$ is
higher than the rank of $\hat{\Cal J}$ by at least $p$, which
shows $\partial \hat{\Cal J} = \hat{\Cal P}$. \qed

\spb We will describe the Frobenius bracket tensor of the
distribution $\Cal P$ in terms of the projective geometry of
$\Cx$. For this, we recall the definitions of the second
fundamental form of a subvariety in the projective space.

  Let $V$ be a complex vector space and $Z \subset \Bbb P V$ be a
subvariety of the projective space. For a smooth point $x \in Z$,
we are going to define a symmetric bilinear form $II_{x,Z}: T_x(Z)
\otimes T_x(Z) \rightarrow N_x(Z;\Bbb P V)$, called the {\it
second fundamental form} of $Z$ at $x$ as follows. Let $\xi:
V\setminus 0 \to \Bbb P V$ be the natural $\Bbb C^*$-bundle. We
denote $\xi^{-1}(x)$ by  $\hat x $ and $\xi^{-1}(Z)$ by $\hat Z$.
 Let $\hat T_x(Z) \subset V$ be the affine tangent space of $Z$
at $x$ which can be naturally identified with the tangent space
$T_{\bar x}(\hat Z)$ of $\hat Z$ at any point $\bar x$ of $\hat
x$. We have natural identifications
$$ \align T_x(Z) &= Hom(\hat x, \hat T_x(Z)/\hat x)\\ N_x(Z; \Bbb
P V) &= Hom(\hat x, V/ \hat T_x(Z)).
\endalign $$
Given two elements $u_o$ and $v_o$ of $T_x(Z)$, choose local
 vector fields
$u$ and $v$ on $Z$ with $u_x = u_o$ and $v_o =v$. Let $\hat u$ and
$\hat v$ be their lifts in an open subset of $\hat Z$. We extend
them to local vector fields in an open subset of $V\setminus 0$
and denote these extended vector fields by the same symbols. In
terms of a linear coordinate system  $x_1, \ldots, x_n$  on $V$,
we can write $$ \align \hat u &= u_1 \frac{\partial}{\partial x_1}
+ \cdots + u_n \frac{\partial}{\partial x_n} \\ \hat v &= v_1
\frac{\partial}{\partial x_1} + \cdots + v_n
\frac{\partial}{\partial x_n}. \endalign$$ We define $$ \align
\hat v(\hat u) &:= \sum_{i,j =1}^{n} v_i \frac{\partial
u_j}{\partial x_i} \frac{\partial}{\partial x_j} \\ \hat u(\hat v)
&:= \sum_{i,j =1}^{n} u_i \frac{\partial v_j}{\partial x_i}
\frac{\partial}{\partial x_j}. \endalign $$ These are local vector
fields on $V \setminus 0$. They are not necessarily tangent to
$\hat Z$, but
$$ \hat v (\hat u) - \hat u (\hat v) = [\hat v, \hat u] $$ is a
lift of the local vector field $[v, u]$ on $Z$. We denote the
value of $\hat v (\hat u)$ at a point $\bar x \in \hat x$ modulo
$\hat T_x (Z) = T_{\bar x}(\hat Z)$ by $\hat{II}_{\bar z}(u_o,
v_o)$. It is easy to see that this vector in $V/ \hat T_x(Z)$ is
independent of the choices of $u, v, \hat u,$ or $\hat v$ and it
depends linearly on the choice of $\bar x \in \hat x$. In other
words, $\bar x \mapsto \hat{II}_{\bar x}(u_o, v_o)$ defines an
element of $Hom( \hat x, V/ \hat T_x(Z)) = N_x(Z; \Bbb P V)$. This
element of $N_x(Z;\Bbb P V)$ is defined to be $II_{x, Z}(u_o,
v_o)$. {}From the fact that $[\hat u, \hat v]$ is tangent to $\hat
Z$, we can see the symmetry $II_{x, Z} (u_o, v_o) = II_{x, Z}
(v_o, u_o)$.

 Another equivalent definition of the second fundamental form is
via the Gauss map.  The {\it Gauss map} $\gamma: Z \dashrightarrow
Gr(\dim Z+1; V)$ is a rational map assigning the tangent spaces to
$Z$ at smooth points of $Z$. In other words, it is defined by
$\gamma(x) = [\hat T_x(Z)]$ at a smooth point $x \in Z$.

 The derivative of the Gauss map $$  d \gamma_x: T_x(Z)
\rightarrow T_{[\hat T_x(Z)]}(Gr(\dim Z+1; V)) = Hom(\hat T_x(Z),
V/\hat T_x(Z)) $$ induces an element of $T_x(Z) \otimes T_x(Z)
\rightarrow N_{x}(Z; \Bbb P V)$ and one can check that our
definition of  $II_{x,Z}$ is just an explicit coordinate
description of this element.
 The following result is classical:

\spb \proclaim{Lemma 2 ([GH, 2.10])} The closures of the fibers of
the Gauss map  are linear subspaces in $V$.
\endproclaim

\spb Now we are ready to describe the Frobenius bracket of $\Cal
P$:

 \spb
\proclaim{Proposition 2} Let
 $\alpha \in \C$. Choose a local complement $\Cal F$ as in
Proposition 1 and  a section $f$ of $\Cal F$ near $\alpha$. Given
two vectors $u_o$ and $ v_o$ in $T_{\alpha}(\Cx)$ with $ x =
\pi(\alpha)$, let $u$ (resp. $ v$) be a local vector field   on a
neighborhood of $\alpha$ in $\C$ tangent to fibers of $\pi$
 such that $u_{\alpha} = u_o$ (resp. $ v_{\alpha} = v_o$).
  Let $[v,f]_{\alpha}$ be the value of the local vector fields
$[v, f]$ on $\C$ at the point $\alpha$. Then the Frobenius bracket
tensor for the distribution $\Cal P$ at $\alpha$
$$ \align [,]:
\bigwedge^2 \Cal P_{\alpha} &\to T_{\alpha}(\C)/\Cal P_{\alpha}\\
&= T_x(X)/\hat T_{\alpha}(\Cx)\\ &= \hat x \otimes N_{\alpha}(\Cx;
\Tx)
\endalign$$ satisfies
$$[u, [ v, f]_{\alpha}] = \pi_* f_{\alpha} \otimes
 II_{\alpha, \Cx}(u, v) \in \hat x
\otimes N_{\alpha}(\Cx; \Tx)$$ where $f_{\alpha} \in \Cal
F_{\alpha}$ denotes the value of $f$ at $\alpha$, which is
contained in $\Cal J_{\alpha}$ so that $\pi_* f_{\alpha} \in \hat
x$. \endproclaim

\spb {\it Proof}. We will work on $\Xi$ as in the proof of
Proposition 1. Choose local coordinate systems $x_{1}, \ldots,
x_{n}$ of $X$ and vertical coordinates $\lambda_{1} = dx_{1},
\ldots, \lambda_{n} = dx_{n}$ of $\Xi$. We choose lifts $\hat v=
\sum_{i}v_{i} \frac{\partial}{\partial \lambda_{i}}$  and $\hat u
= \sum_{i} u_{i} \frac{\partial}{\partial \lambda_{i}}$. As in the
proof of Proposition 1, we may choose $\hat f = \sum_{i} f_{i}
\frac{\partial}{\partial \lambda_{i}} + \sum_{j} \lambda_{j}
\frac{\partial}{\partial x_{j}}$. Then $$ \align [\hat v, \hat f]
&= \sum_{i} v_{i} \frac{\partial}{\partial x_{i}} +
\sum_i (\hat v (f_i)- \hat f(v_i)) \frac{\partial}{\partial \lambda_i} \\
[\hat u, [\hat v, \hat f]] &\equiv \sum_{i} \hat u (v_{i})
\frac{\partial}{\partial x_{i}} \mod \hat{\Cal V}. \endalign $$
Restricting $\hat{u}$ and $\hat{v}$ to $\hat{\C}_x$, we wee that
$[u, [v,f]] = \hat{II}_{\alpha, \Cx}(u, v)$. \qed

\spb Using Lemma 2, we can define a sub-distribution $\Cal G$ of
$\Cal V$ which consists of tangent vectors to fibers of the Gauss
map for $\Cx$ as $x$ varies. This distribution is integrable. An
immediate consequence of Proposition 2 is

\spb \proclaim{Proposition 3} In the notation above, $\Cal G =
Ch(\Cal P) \cap \Cal V$. \endproclaim

\spb This enables us to describe $Ch(\Cal P)$ in terms of the
Gauss map for $\Cx$ as follows.

 \spb \proclaim{Proposition 4}
If there exists a complement $\F$ to $\Cal V$ in $\Cal J$ in an
open subset of $\C$ so that $\Cal F \subset Ch(\Cal P)$, then
$Ch(\Cal P) = \F + \Cal G + [\F, \Cal G]$ on that open subset. In
this case, if the rank of ${\Cal G}$ is $k-1$ for some $k >0$,
then the rank of $Ch({\Cal P})$ is $2k-1$.
\endproclaim

\spb {\it Proof}. {}From Proposition 3 and the fact that $Ch(\Cal
P)$ is closed under Lie bracket, the inclusion $\F + \Cal G + [\F,
\Cal G] \subset Ch(\Cal P)$ is immediate.
 We know that any local section of $Ch(\Cal P)$ can be written as
$f_1 + v_1 + [f_2, v_2]$ for some local sections $f_1, f_2$ of
$\F$ and $v_1, v_2$ of $\Cal V$ from Proposition 1. We want to
show that this local section is in $\F + \Cal G + [\F, \Cal G]$.
It suffices to show that $v_2$ is a section of $Ch(\Cal P)$. In
fact, if $v_2$ is a section of $Ch(\Cal P)$,   it is a section of
$\Cal G$ from Proposition 3 and so $[f_2, v_2] \in [\Cal F, \Cal
G]$. As a consequence, we have $v_1 \in Ch(\Cal P) \cap \Cal V =
\Cal G$, which proves Proposition 4. To show that $v_2$ is a
section of $Ch(\Cal P)$, we need to check that $[v_2, h]$ is a
section of $\Cal P$ for any local section $h$ of $\Cal P$. {}From
Proposition 1, we can set $h= f_3 + v_3 + [f_4, v_4]$. Then
$$\align [v_2, h] &= [v_2, f_3 + v_3+ [f_3, v_4]
\\ &= [v_2, f_3] + [v_2, v_3] + [ v_2, [f_3, v_4]]\\ &= [v_2, f_3]
+[v_2, v_3] + [v_4, [f_3, v_2]] + [f_3, [v_2, v_4]]. \endalign $$
So it suffices to show that $[v_4, [f_3, v_2]]$ is a section of
$\Cal P$. Since $f_1 + v_1 + [f_2, v_2]$ is a section of $Ch(\Cal
P)$, we see that  $[v_4, f_1 + v_1 + [f_2, v_2]]$ is  a section of
$\Cal P$. Thus $[v_4, [f_2, v_2]]$ is a section of $\Cal P$. This
implies that $[v_4, [f_3, v_2]]$ is a section of $\Cal P$ because
$f_3 = \zeta f_2$ for some local holomorphic function $\zeta$.

Finally, the statement about the ranks follows from the local
coordinate expression of $[\F, \Cal G]$ as in the proof of
Proposition 1. \qed

\spb \proclaim{ Proposition 5} Suppose there exists a local
complement $\F$ to $\Cal V$ in $\Cal J$ as in Proposition 4. Let
$S$ be a local leaf of $Ch(\Cal P)$ in an open neighborhood of the
domain of definition of $Ch(\Cal P)$ where $\F$ is a well-defined
foliation.  Then $S$ is an open subset in $\Bbb P T(\pi(S))$.
\endproclaim

\spb {\it Proof}. Since the leaves of $\F$ in $S$ are sent to
holomorphic curves in $\pi(S)$, we see that $S \subset \Bbb P
T(\pi(S))$. Let $k-1$ be the rank of $\Cal G$ so that $\dim S =
2k-1$. Note that the intersection of  $S$ with a fiber
$\pi^{-1}(x)$ is an open subset of the projective space in $\Cx$
which is a fiber of the Gauss map of $\Cx$. Thus we get $\dim(
\pi(S)) =k$, which implies that $\dim(S) = \dim \Bbb P T(\pi(S))$.
\qed

\mpb
 \noindent {\bf \S 3. Birationality of the tangent map}

\spb Let $X$ be a uniruled projective manifold and $\K$ be a
minimal component. Let $\rho: \U \rightarrow \K$ and $\mu: \U
\rightarrow X$ be the associated universal family morphisms. Since
$\tau: \U \dashrightarrow \Cal C$ is generically finite, there
exists a Zariski open subset $\U^o$ such that $\tau|_{\U^o}$
   is an etale morphism. The fibers of $\rho: \U \rightarrow \K$
can be   regarded as a foliation by curves on an etale cover of an
open subset of $\C$. This will be called the {\it tautological
foliation} on $\C$ and denoted by $\F$. This is a multi-valued
foliation on $\C$.

\spb \proclaim{Proposition 6} The tautological foliation $\F$ is a
univalent foliation on $\C$ if and only if the tangent map $\tau:
\U \rightarrow \C$ is birational. \endproclaim

\spb {\it Proof}. The univalence of $\F$ when $\tau $ is
birational is obvious. If $\F$ is univalent, its leaves are curves
on $\C$ whose images in $X$ under the projection $\pi: \C
\rightarrow X$ are just members of $\K$. Thus $\tau$ must be
birational from the generic injectivity of the natural map $\K
\rightarrow \text{Chow}(X)$. \qed

\spb Now let us apply Section 2 to the subvariety $\C$ of $\Bbb P
T(X)$. We have natural foliations $\Cal V$, $\Cal T$ and $\Cal P$
on $\C$ which are completely determined by the inclusion $\C
\subset \Bbb P T(X)$.

\spb \proclaim{Proposition 7} In a neighborhood of a general point
of $\C$, choose a complex-analytic open subset of $\C$ so that a
univalent choice of the values of  the tautological foliation $\F$
can be made on that open subset. Denote this univalent foliation
by the same symbol $\F$. Then $\Cal V + \F = \Cal T$ in that open
subset.
\endproclaim

\spb {\it Proof}. By the definition of $\Cal T$, $\Cal T_\alpha =
(d \pi)^{-1} (\Bbb C \alpha)$ at $\alpha \in \C$, it is obvious
that $\F \subset \Cal T$. Thus $\Cal V + \F = \Cal T$ follows from
$\F \nsubseteq \Cal V$ and $ rank(\Cal V) + rank(\F) = rank(\Cal
T).$ \qed

\spb On $\K$, we have a natural distribution $\Cal R$ defined as
follows. At a general point $[h]$ of $\K$ corresponding to a
rational curve $h: \P \rightarrow X$, the bundle $h^*T(X)$ on $\P$
splits as $\Cal O(2) \oplus \Cal O(1)^p \oplus \Cal O^{n-1-p}$ for
some non-negative integer $p$.  Members of $\K$ having such
splitting type of $T(X)$ are called {\it standard rational
curves}. By elementary deformation theory, $p$ is the fiber
dimension of $\mu: \U \rightarrow X$ and the tangent space to $\K$
is
$$T_{[h]}(\K) = H^0(\P, h^* T(X)) / H^0 (\P, T(\P)) \cong H^0( \P, \O(1)^{p}
\oplus \O^{n-1-p}).$$ The subspace of $T_{[h]}(\K)$ corresponding
to the subspace $H^0(\P, \O(1)^p)$ is determined independent of
the choice of the isomorphism $h^*T(X) \cong \O(2) \oplus \O(1)^p
\oplus \O^{n-1-p}$. This subspace will be defined to be $\Cal
R_{[h]}$. Thus $\Cal R$ is a distribution of rank $2p$ on $\K$
whose domain of definition includes the open subset consisting of
standard rational curves belonging to $\K$.

\spb \proclaim{Proposition 8} The pull-back of $\Cal R$ by $\rho$
agrees with the pull-back of $\Cal P$ by $\tau$: $d \tau^* \Cal P
= d \rho^* \Cal R.$ In particular, $\F \subset Ch(\Cal P)$.
\endproclaim

\spb {\it Proof}. At a general point $\alpha \in \C$ corresponding
to the tangent direction to a general standard minimal rational
curve $h: \P \rightarrow X$ with $h(o)=x$, we have the natural
identifications
$$ T_{\alpha}(\C) = H^0(\P, h^*T(X))/H^0(\P, T(\P) \otimes \bold m_o)$$
$$ T_x(X) = H^0(\P, h^*T(X))/H^0(\P, h^*T(X) \otimes \bold m_o).$$
Under these identifications,  the projection $d \pi: T_{\alpha}
\rightarrow T_x(X)$ corresponds to taking the value of a section
in $H^0(\P, h^*T(X))$ at the point $x$.  By definition, $\Cal
P_\alpha = (d \pi)^{-1}(\hat{T}_{\alpha}(\Cx))$ at a generic point
$\alpha \in \C$. By  elementary deformation theory,
$\hat{T}_{\alpha}(\Cx)$ is naturally isomorphic to $(\O(2) \oplus
\O(1)^p)$-part of the splitting $h^*T(X) \cong \O(2) \oplus
\O(1)^p \oplus \O^{n-1-p}$. Thus $\Cal P_{\alpha}$ corresponds to
the sections in $H^0(\P, h^*T(X))$ whose values at $x$ lie in the
$\O(2)\oplus \O(1)^p$-part of the splitting. Thus it is exactly
the lift of $\Cal R_{[h]}$ by $\rho$. \qed

\spb \proclaim{Proposition 9} Let $S$ be a general leaf of
$Ch(\Cal P)$. Then $\pi(S)$ is a quasi-projective variety and
contains a smooth Zariski open subset $\V \subset \pi(S)$ such
that $\Bbb P T(\V)$ is a Zariski dense open subset of $S$. Let
$\K^S$ be the subvariety of $\K$ consisting of members of $\K$
lying on the closure $\bar{S}$ of $\pi(S)$. Then $\K^S$ is a
minimal component of the irreducible projective variety $\bar{S}$.
Moreover, the corresponding total space of minimal rational
tangents $\C^S$ agrees with $\Bbb P T(\bar{S})$.
\endproclaim

\spb {\it Proof}. By Proposition 7 and Proposition 8, we can apply
Proposition 5 here. The result follows from the fact that the
leaves of $\F$ correspond to minimal rational curves. \qed

\spb Now we recall the following result of
Cho-Miyaoka-Shepherd-Barron. Note that  minimal components and
their total spaces of minimal rational tangents are defined for
any irreducible projective variety $X$ in Section 1.

\spb \proclaim{Proposition 10 ([CMS])} Let $X$ be an irreducible
normal projective variety of dimension $n$. Suppose there exists a
minimal component $\K$ on $X$ such that the total space of minimal
rational tangents agrees with $\Bbb PT(X)$. Then there exists a
finite morphism $ \Bbb P^n \rightarrow X$ which is etale over
$X-Sing(X)$ such that the members of $\K$ are just the images of
lines in $\Bbb P^n$. In particular, $X \cong \Bbb P^n$ if $X$ is
smooth.
\endproclaim

\spb A simple consequence is the following special case of Theorem
1.

\spb \proclaim{Proposition 11} Let $X$ be an irreducible normal
projective variety which has a minimal component $\K$ such that
the total space of minimal rational tangents $\C$ agrees  with
$\Bbb PT(X)$. Then

(i) the tangent map $\tau: \U \rightarrow \Cal C$ is birational
and

(ii) $\K$ is the only minimal component of $X$.
\endproclaim

\spb We are ready to prove Theorem 1 and 2.

 \spb {\it Proof of Theorem 1 and Theorem 2}. Let us start with Theorem 1.
 By Proposition 6, it
suffices to show that the tautological foliation $\F$ is
univalent. By Proposition 8, $\F$ is contained in $Ch(\Cal P)$ and
it suffices to show that it is univalent on a general leaf $S$ of
$Ch(\Cal P)$. But by Proposition 9, $\F$ restricted to $S$ is the
tautological foliation for the minimal component $\K^{S}$ of
$\bar{S}$ with the total space of minimal rational tangents $\C^S$
equal to  $\Bbb P T(\bar{S})$. Thus $\F$ is univalent on $S$ by
Theorem Proposition 10 (i). Theorem 2 follows from the fact that
$\K^{S}$ is the only minimal component of  $\bar{S}$, by
Proposition 10 (ii). \qed

\spb We will finish this section with a few results about $\pi(S)$
in the notation of Proposition 9. A subvariety $Y$ in $X$ is
called a {\it Cauchy subvariety} of the minimal component $\K$ if
it is the closure of $\pi(S)$ for a general leaf $S$ of $Ch(\Cal
P)$ in the notation of Proposition 9. The following is an
immediate consequence of Proposition 10.

\spb \proclaim{Proposition 12} Let $X$ be a uniruled projective
manifold and $Y \subset X$ be a Cauchy subvariety with respect to
a choice of a minimal component $\K$. Let $\tilde{Y}$ be the
normalization of a Cauchy subvariety $Y$. Then there exists a
finite morphism $ \Bbb P^d \rightarrow \tilde{Y}$ which is etale
over $\tilde{Y}-Sing(\tilde{Y})$ such that the members of $\K$
lying on $Y$ are just the images of lines in $\Bbb P^d$.
\endproclaim

\spb Corollary 1 has the following consequence.

 \spb \proclaim{Proposition 13} Let $X$ be a uniruled projective
manifold and $\K$ be a minimal component. Suppose the variety of
minimal rational tangents $\Cx$ is non-linear. Then for each
component $\Cx^1$ of $\Cx$, the intersection of the closures of
the fibers of the Gauss map for $\Cx^1$ is empty.
\endproclaim

\spb {\it Proof}. Suppose there exists a point $\alpha \in \Cx^1$
which is the intersection of the closures of the fibers of the
Gauss map. The normalization $\K^1_x$ of $\Cal C_x^1$ is smooth by
Corollary 1. The fibers of Gauss map give subvarieties of the
smooth projective variety $K^1_x$ passing through a point
$\hat{\alpha} \in \K_x^1$ over $\alpha$. The pull-back of $\Cal
O(1)$ bundle on $\Bbb PT_x(X)$ to $\K_x^1$ gives an ample line
bundle $L$ with respect to which the fibers of the Gauss map are
linear subspaces. It follows that there exists a family of
rational curves through $\hat{\alpha}$ covering the whole smooth
variety $\K_x^1$ such that each members have degree 1 with respect
to the ample line bundle $L$. By Proposition 10, $\K_x^1$ must be
a projective space and then $\Cal C^1_x$ must be linear. \qed

\spb  It is well-known that a positive-dimensional family of
minimal rational curves passing through a general point $x \in X$
do not have a common intersection other than $x$
(`bend-and-break'). A weaker form of this fact can be generalized
to Cauchy subvarieties as follows.

\spb \proclaim{Proposition 14} Let $X$ be a uniruled projective
manifold and $\K$ be a minimal component. Assume that $\Cal C_x$
is non-linear for a general point $x \in X$. Let $\Cal C_x^1 $ be
a component of $\Cal C_x$. Then the irreducible family of the
Cauchy subvarieties passing through $x$ whose tangent spaces lie
in $\Cal C_x^1$ do not have a common intersection point other than
$x$.
\endproclaim

\spb {\it Proof}. Note that the non-linearity of $\C_x$ implies,
by Lemma 2, that there are a positive-dimensional family
$\{Y_s\}$, where $s$ is a parameter, of Cauchy subvarieties
passing through $x$ whose tangent spaces lie in $\Cal C_x^1$.
Suppose there exists a common intersection point $y$ different
from $x$. By Proposition 12, for each $s$, there exists a member
$C_s$ of $\K$ passing through $x$ and $y$ lying on the Cauchy
subvariety $Y_s$. Since there cannot exist a positive-dimensional
family of minimal rational curves passing through $x$ and $y$,
$C_s$ cannot be a family of distinct members of $\K$. It follows
that a single member $C=C_s$ belongs to each $Y_s$ and  $\alpha:=
\Bbb P T_x(C_s)$ is in the intersection of the fibers of the Gauss
map of $\Cal C_x^1$, a contradiction to Proposition 13. \qed

\mpb
 \noindent {\bf \S 4. Rigidity of generically finite morphisms to
Fano manifolds with non-linear varieties of minimal rational
tangents}

\spb A weaker form of Theorem 3 was already proved in [HM2] where
the stronger assumption that $\Cx$ has generically finite Gauss
map was needed. The proof of Theorem 2 is a generalization of the
argument in [HM2] using the result of Section 3. We start with
recalling two results from [HM2].

 \spb \proclaim{Lemma 3 ([HM2, Lemma 4.2])} Let $\pi:Y
\rightarrow X$ be a generically finite morphism from a normal
variety $Y$ onto a Fano manifold $X$ with Picard number 1. Suppose
for a  general member $C \subset X$ belonging to a chosen minimal
 component ${\Cal K}$, each component of the inverse image
$\pi^{-1} (C)$ is birational to $C$ by $\pi$. Then $\pi:Y
\rightarrow X$ itself is birational.
\endproclaim

 \spb \proclaim{Proposition 15 ([HM2, Proof of Theorem 1.4])}
  Let $\{ X_t, t \in \Delta \}$ be
a regular family of uniruled projective manifolds. Given a minimal
component $\K_0$ of $X_0$ with the total space of variety of
minimal rational tangents $\Cal C_0 \subset \Bbb P T(X_0)$, there
exists $\epsilon>0$ and a family of minimal components $\K_t$ of
$X_t, 0<|t| < \epsilon$, with the total space  of minimal rational
tangents $\Cal C_t \subset \Bbb P T(X_t)$ having the following
property: given a family of generically finite morphisms $f_t:X'
\rightarrow X_t$ from a fixed projective manifold $X'$ and a
general point $x\in X'$, the family of subvarieties
$(df_t)_x^{-1}(\Cal C_{t, f_t(x)})$ in $\Bbb P T_x(X')$ is a
constant family, i.e.
$$(df_t)^{-1}(\Cal C_{t, f_t(x)}) = (df_0)^{-1}(\Cal C_{0,
f_0(x)})$$ for all $|t| < \epsilon$.
\endproclaim

\spb The condition of the non-linearity of the variety of minimal
rational tangents is used in the following manner.

 \spb \proclaim{ Proposition 16} Let $X$ be a Fano manifold of
Picard number 1. Suppose there exists a minimal component ${\Cal
K}$  such that the variety of minimal rational tangents ${\Cal
C}_{x}$ at a general point $x \in X$ is non-linear of dimension
$p>0$. Let $X'$ be another Fano manifold of Picard number 1 with a
minimal dominating component $ {\Cal K}'$. Assume that the variety
of minimal rational tangents ${\Cal C}'_{x'} \subset {\Bbb
P}T_{x'}(X')$ at a general point $x'$ has dimension $p$.
  If there exists a quasi-projective variety $U$ and  etale morphisms
  $e:U \rightarrow
X$ and  $\varphi:U \rightarrow X'$ preserving varieties of minimal
rational tangents in the sense that $\varphi_{*}e^{*}({\Cal
C}_{e(y)})  \subset {\Cal C}'_{\varphi (y)}$ for all $y \in U$,
then there exists a birational map $\Phi:X \rightarrow X'$ such
that $\varphi= \Phi \circ e$. \endproclaim

\spb {\it Proof}.  We define an equivalence relation on $U$ by
$y_{1} \sim y_{2}$ if $e(y_{1}) = e(y_{2})$ and $\varphi(y_{1}) =
\varphi(y_{2})$. In the quotient space $U/\!\!\!\sim$, we can find
a Zariski-open subset $\tilde{U}$ so that the induced morphisms
$\tilde{e}:\tilde{U} \rightarrow X$ and $\tilde{\varphi} :
\tilde{U} \rightarrow X'$ are \'etale, preserving varieties of
minimal rational tangents in the above sense. For a general Cauchy
subvariety $Y$ of  ${\Cal K}$, any component $Y_{1}$ of
$\tilde{e}^{-1}(Y)$ is sent to a Cauchy subvariety of $\Cal K'$
because Cauchy subvarieties are completely determined by the
geometry of the total space of  minimal rational tangents in
${\Bbb P}T(X)$ and $\Bbb P T(X')$.

We claim that $\tilde{e}$ is injective. Suppose not. Then by Lemma
3, for a general member $C$ of ${\Cal K}$, there exists a
component $\hat{C}$ of $\tilde{e}^{-1}(C)$ so that
$\tilde{e}|_{\hat{C}}: \hat{C} \rightarrow C$ is not birational.
Choose a general point $x \in C$ and let $y_{1} \neq y_{2}$ be
points of $\hat{C}$ over $x$. By the construction of $\hat{U}$ we
can assume that $\hat{\varphi}(y_{1}) \neq \hat{\varphi}(y_{2})$.
 Since $C$ is general, we can assume that for
 $p$-dimensional deformations $C_{t}$ of $C$ fixing $x$,  the collection
of some  components $\hat{C}_{t}$ of $\tilde{e}^{-1}(C_{t})$ gives
a $p$-dimensional deformations of $\hat{C}$ fixing $y_{1}$ and
$y_{2}$. For the Cauchy subvariety $S_t$ containing $C_t$,
$\tilde{e}^{-1}(S_t)$ contains an irreducible component $Z_t$
containing both $y_1$ and $y_2$. Thus we have a component ${\Cal
C}^1$ of the variety of minimal rational tangents ${\Cal
C}'_{\hat{\varphi}(y_{1})}$ and a component $\Cal C^2$ of the
variety of minimal rational tangents ${\Cal
C}'_{\hat{\varphi}(y_{2})}$ such that the  Cauchy subvarieties
corresponding to $\Cal C^1$ and the Cauchy subvarieties
corresponding to $\Cal C^2$ are equal.  This is a contradiction to
Proposition 14. \qed

\spb{\it Proof of Theorem 3}. Suppose that $f$ is  birational.
Then $f$ is biholomorphic over ${\Cal X} -{\Cal Z}$ where ${\Cal
Z}$ is
 a subvariety of codimension $\geq 2$.
On ${\Cal X}'$ we have the vector field ${\Cal V}$ lifting
 $\frac{d}{d t}$ on $\Delta$. Its push-forward $f_{*}{\Cal V}$ is a
vector field on ${\Cal X} -{\Cal Z}$. By Hartogs, we can extend it
to a vector field  on ${\Cal X}$ which generates the required
family of biholomorphic morphisms.

Suppose that $f_{t}$ is not birational, but generically $d$-to-1.
We will construct a new projective manifold $\hat{X}$, a
generically finite dominant rational map
 $\nu:X' \rightarrow \hat{X}$ of degree $d$ and a holomorphic family of
  birational morphisms
$g_{t}: \hat{X}  \rightarrow X_{t} $ for small $t$, so that $f_{t}
= g_{t} \circ \nu$ over general points of $X'$.
 Then the proof is reduced to the  birational case, which is already settled above.

Let $\Cal K^0$ be the minimal component on $X_0$ with the total
space  of minimal rational tangents $\Cal C^0 \subset \Bbb
PT(X_0)$ whose fiber over a general point of $X_0$ consists of
non-linear subvarietes.  Applying Proposition 15, there exists
some $\epsilon
>0$ and a minimal component $\Cal K_t$ of $X_t$ for $0 < |t| <
\epsilon$ such that  the total space of minimal rational tangents
$\Cal C^t$ satisfies  $$ (df_t)^{-1}_y ( \Cal C^t_{f_t(y)}) =
(df_0)^{-1}_y ( \Cal C^0_{f_0(y)}).$$

 Fix such a small $t=s$. There exists an
open subset $U \subset X'$ which is etale over $X_0$ and $X_s$ by
$f_0$ and $f_s$. The above equality of the inverse image of the
varieties of minimal rational tangents at $f_s(y)$ and $f_0(y)$
implies that the hypothesis of Proposition 16 is satisfied for the
etale morphisms $f_s|_U$ and $f_0|_U$. Thus there exists a
birational map $\Phi: X_0 \dashrightarrow X_s$ such that $f_s|_U =
\Phi \circ f_0|U$.

  We say that a reduced 0-cycle $y_{1}+
\cdots +y_{d}$ of length $d$
 on $X'$ is a special cycle
 if  $f_0(y_{1}) =
\cdots = f_0(y_{d})$ and $f_{s}(y_{1}) = \cdots = f_{s}(y_{d})$.
{}From the existence of the birational map $\Phi$, general fibers
of $f_0$ and $f_s$ are special cycles. The set of all special
cycles  on $X'$ gives a constructible subset of the Hilbert scheme
of 0-dimensional subschemes of $X'$. We  can find an irreducible
component
 of this set so that
the corresponding cycles cover an open set of $X'$. Let $B$ be the
closure of that component and let
 $\sigma: A \rightarrow B$ and $ \lambda: A \rightarrow X'$ be
 the universal family morphisms so that
$\sigma$ is flat and of degree $d$. We claim that $\lambda$ is
birational. Otherwise we have two distinct
 special cycles of degree $d$ containing
a given general point $y$ of $X'$, which is absurd because $f_0$
and $f_s$ are of degree $d$.

Let $\hat{X}$ be a desingularization of $S$ and
$\hat{\sigma}:\hat{A} \rightarrow \hat{X}, \hat{\lambda} : \hat{A}
\rightarrow X'$ be the induced morphisms. Define the rational map
$\nu: X' \dashrightarrow \hat{X}$ as $\nu := \hat{\sigma}
 \circ \hat{\lambda}^{-1}$. Then $\nu$ is a generically finite
dominant rational map of degree $d$. Consider the morphism
$\hat{f}_{s} = f_{s} \circ \hat{\lambda}$ from $\hat{A}$ to
$X_{s}$. {}From the definition of special cycles,
 a general fiber of
$\hat{\sigma}$ is contained in a fiber of $\hat{f}_{s}$. Thus each
fiber of $\hat{\sigma}$ is contained in a fiber of $\hat{f}_{s}$
 by the flatness of $\hat{\sigma}$. It follows that
 we get a birational morphism $g_{s}: \hat{X} \rightarrow X_{s}$
 satisfying $f_{s} = g_{s} \circ \nu$. Since $ X', X_{s}$ are all projective,
 it is easy to see that $\{g_{s}\}$ is  a holomorphic family. \qed

\spb Let $X$ be a Fano manifold of Picard number 1.  Suppose $X'$
is a projective manifold with non-zero holomorphic vector fields.
Given a generically finite morphism $f:X' \rightarrow X$, we get a
deformation $f_t:X' \rightarrow X$ of $f$ obtained by composing it
with the 1-parameter subgroup of automorphisms of $X'$ generated
by the holomorphic vector fields. In this situation,  Theorem 3,
with $\Cal X = X \times \Delta$, implies the following.

 \spb
\proclaim{Corollary 2} Let $X$ be a Fano manifold of Picard number
1 which has a minimal  component with non-linear variety of
minimal rational tangents. Then for any projective variety $X'$
and a surjective generically finite morphism $f: X' \rightarrow
X$, any holomorphic vector field on $X'$ descends to a holomorphic
vector field on  $X$ such that $f$ is equivariant with respect to
the 1-parameter groups of automorphisms of $X'$ and $X$ generated
by the holomorphic vector fields.
\endproclaim

\spb {\it Proof of Theorem 4}. Since the Picard number of $G/P$ is
1, $f$ is finite and $X$ is a Fano manifold of Picard number 1.
Assume that $X \not\cong \Bbb P_n$ and pick a minimal component
$\Cal K$ on $X$ with $\C \neq \Bbb PT(X)$ by Proposition 10.

First assume that $\Cx$ is non-linear at general $x \in X$. By
Corollary 2, vector fields on $G/P$ descend to vector fields on
$X$ by $f$. If $R \subset G/P$ is a non-empty ramification divisor
of $f$ and $B=f(R)$, then the image of any tangent vector of $G/P$
at a point $x$ on $R$ is tangent to $B$ whenever $f(x)$ is a
smooth point of $B$. Thus all vector fields of $G/P$ are sent to
vector fields on $X$ tangent to $B$ at smooth points of $B$. This
implies that integral curves of any vector field of $G/P$ through
a point of $R$ are sent to curves in $B$. In other words, the
integral curves through a point of $R$ remain in $R$, a
contradiction to the homogeneity of $G/P$.

Now we may assume that $\Cx$ is linear. Then a component of
$(df)^{-1}(\Cal C_{f(s)}) \subset \Bbb P T_s(S)$ at a general
point $s \in S$ is a linear space invariant under the isotropy
group $P$ by [HM1, Section 1], defining a $G$-invariant
distribution on $G/P$. By translating an irreducible component of
the inverse image of a Cauchy subvariety, we get a proper
$G$-invariant foliation on $G/P$ with compact leaves, a
contradiction to the assumption that $G/P$ has Picard number 1.
\qed

 \mpb {\bf References}

\spb

[CMS] Cho, K., Miyaoka, Y. and Shepherd-Barron, N. I.:
Characterization of projective space and applications to complex
symplectic manifolds. in {\it Higher dimensional birational
geometry (Kyoto, 1997)},  Adv. Stud. Pure Math. {\bf 35}  (2002)
1-88

[GH] Griffiths, P. and Harris, J.: Algebraic geometry and local
differential geometry. Ann. scient. \'Ec. Norm. Sup. {\bf 12}
(1979) 355-432

[HM1] Hwang, J.-M. and Mok, N.: Holomorphic maps from rational
homogeneous spaces of Picard number 1 onto projective manifolds.
Invent. math. {\bf 136} (1999) 209-231

[HM2] Hwang, J.-M. and Mok, N.: Cartan-Fubini type extension of
holomorphic maps for Fano manifolds of Picard number 1.  Journal
Math. Pures Appl. {\bf 80} (2001) 563-575

[Ke] Kebekus, S.: Families of singular rational curves.
 J. Algebraic Geom. {\bf 11} (2002) 245--256

[KK] Kebekus, S. and Kovacs, S.:Are minimal degree rational curves
determined by their tangent vectors? preprint alg-geom/0206193

 [Ko] Koll\'ar, J.: {\it Rational curves on algebraic
varieties.} Ergebnisse der Mathematik und ihrer Grenzgebiete, 3
Folge, Band 32, Springer Verlag, 1996

\mpb\noindent Jun-Muk Hwang (jmhwang\@ns.kias.re.kr)
\newline
Korea Institute for Advanced Study, 207-43 Cheongryangri-dong
\newline
Seoul 130-722, Korea

\mpb\noindent Ngaiming Mok (nmok\@hkucc.hku.hk)
\newline
Department of Mathematics, The University of Hong Kong \newline
Pokfulam Road, Hong Kong

 \bye